\documentclass[11pt]{article}

\usepackage[margin=1in]{geometry}
\usepackage{amsmath,amsfonts,amsthm,latexsym,xspace}
\usepackage{graphicx}
\usepackage{dsfont}
\usepackage{booktabs}
\newtheorem{thm}{Theorem}[section]

\newtheorem{lem}[thm]{Lemma}

\theoremstyle{definition}

\theoremstyle{remark}
\newtheorem{rem}[thm]{Remark}

\newtheoremstyle{quote}{2pt}{-2pt}{\it}{}{\bf}{\mbox{\rm .}}{0.5em}{}%
\theoremstyle{quote}

\newcommand{\eps}{\varepsilon}
\newcommand{\M}{\mathcal{M}}
\newcommand{\E}{\mathbb{E}}
\newcommand{\CH}{\mathcal{H}}
\newcommand{\C}{\mathcal{C}}
\newcommand{\CV}{\mathcal{V}}
\newcommand{\CE}{\mathcal{E}}

\newcommand{\dd}{\mathrm{d}}

\title{{\bf Metric Construction, Stopping Times\\
  and Path Coupling}}
\author{Magnus Bordewich\thanks{School of
Computing, University of Leeds, Leeds LS2 9JT, UK. Email:
\texttt{\{dyer,magnusb\}@comp.leeds.ac.uk}.},\ \, Martin Dyer${}^*$ and Marek
Karpinski\thanks{Dept. of Computer Science, University of Bonn, 53117 Bonn,
Germany. Email: \texttt{marek@cs.uni-bonn.de}.}}
\date{November 15, 2005}
\begin{document}
\begin{titlepage}
\maketitle

\begin{abstract}
In this paper we examine the importance of the choice of metric in
path coupling, and the relationship of this to \emph{stopping time
analysis}. We give strong evidence that stopping time analysis is no
more powerful than standard path coupling. In particular, we prove a
stronger theorem for path coupling with stopping times, using a
metric which allows us to restrict analysis to standard one-step
path coupling. This approach provides insight for the design of
non-standard metrics giving improvements in the analysis of specific
problems. 

We give illustrative applications to hypergraph
independent sets and SAT instances, hypergraph colourings and
colourings of bipartite graphs. In particular we prove rapid mixing for Glauber dynamics on independent sets in hypergraphs whenever the minimum edge size $m$ and degree $\Delta$ satisfy $m\geq \Delta + 2$, and for all edge sizes when $\Delta=3$. Previously rapid mixing was only known for $m\geq 2\Delta+1$. This result leads to approximation schemes for monotone SAT formulae in which the maximum number of occurrences of a variable ($\Delta$) and the minimum number of variables per clause ($m$) satisfy the same condition. For Glauber dynamics on proper colourings of 3-uniform hypergraphs we prove rapid mixing whenever the number of colours $q$ is at least $\big\lceil
\tfrac{3}{2}\Delta+1 \big\rceil$. Previously the best known result was for $q\geq 1.65 \Delta$ and $\Delta\geq \Delta_0$ for some large $\Delta_0$. Finally we prove rapid mixing of scan dynamics (where the order of vertex updates is deterministic) for proper colourings of bipartite graphs whenever $q>f(\Delta)$, where $f(\Delta)\rightarrow \beta \Delta$, as $\Delta\rightarrow \infty$, and $\beta$ satisfies $\tfrac{1}{\beta}e^{1/\beta}=1,\ \ (\beta\approx 1.76)$. This gives rapid mixing with fewer colours than Vigoda's $11\Delta/6$ bound~\cite{V99}, whenever $\Delta \geq 31$, and equals this bound for $\Delta\geq 14$.
\end{abstract}
\end{titlepage}
\section{Introduction}\label{sec:intro}

Path coupling~\cite{BD97b} has proved to be a useful technique for
analysing Markov chains. Analysis is carried out relative to a
chosen \emph{metric} on the state space, for example the Hamming
distance on the independent sets in a graph or hypergraph. The
limitations of the analysis are always caused by certain ``bad''
configurations. But these configurations may be unlikely in a typical
realisation of the chain. Consequently, path coupling has been
augmented by other techniques, such as \emph{stopping time}
analysis. See \cite{BDK05a,DGGJM01,HV04,MN02} for some applications
of this technique. A general theorem for applying stopping times was
proved in \cite{HV04}, and improved somewhat in \cite{BDK05a}.

The stopping time approach is applicable when the bad configurations
have a reasonable probability of becoming less bad as time passes.
For example,  the bad configurations for the Glauber dynamics on
hypergraph independent sets involve almost full edges containing the
change vertex. (See~\cite{BDK05a} for details.) However, it seems
likely that the number of occupied vertices in these edges will have
been reduced before we must either increase or decrease the distance
between the coupled chains. This observation allows a greatly
improved analysis~\cite{BDK05a}.

The stopping time approach is a multistep analysis, and appears to
give a powerful extension of path coupling. However, in this paper
we provide  strong evidence that the stopping time approach is no
more powerful than single-step path coupling. We observe that, in
cases where stopping times can be employed to advantage, equally
good or better results can be achieved by using a suitably tailored
\emph{metric} in the one-step analysis. The intuition behind the
choice of metric will be illustrated with several examples.

In fact, our first example is a proof of a theorem  for path
coupling using stopping times, relying on a particular choice of
metric which enables us to work with the standard one-step path
coupling. The resulting theorem is stronger than those in
\cite{BDK05a,HV04}.  The proof implies that all results obtained
using stopping times can just as well be obtained using standard
path coupling and the right choice of metric. This does not
immediately imply that we can abandon the analysis of stopping
times. Determining the metric used in our proof involves bounding
the expected distance at a stopping time. However the proof does
suggest that it may be better to carry out one-step analysis using a
metric indicated directly by the stopping time intuition.

With this insight, we revisit the Glauber dynamics for hypergraph
independent sets (or equivalently, satisfying assignments of
monotone SAT formulas), and hypergraph colourings, analysed
in~\cite{BDK05a} using stopping times. We find that we are able to
obtain considerably stronger results than those obtained
in~\cite{BDK05a}, using metrics inspired by the stopping times
considerations but then optimised to give the best results. The
technical advantage arises from the possibility of using
linearity of expectation where stopping time analysis must use
concentration inequalities and union bounds.

We note that this paper does not contain the first uses of
``clever'' metrics with path coupling. See \cite{BD98,LV99} for
examples. But we do give the first widely applicable rationale for
choosing a good metric. While there have been instances in the
literature of optimising the \emph{chain}~\cite{DG00a,V99}, the only
previous analysis of which we are aware which uses optimisation of
the \emph{metric} appeared in~\cite{LV99}.

The organisation of the paper is as follows. In section
\ref{sec:stopping} we  prove a better stopping time theorem than
previously known, using only standard path coupling. In section
\ref{sec:indsets} we give our improved results for sampling
independent sets in hypergraphs, and in section \ref{sec:sat},
applications to counting the number of satisfying assignments in
monotone SAT formulas. In section \ref{sec:3uniform} we give
improved results for sampling colourings of $3$-uniform hypergraphs.
Finally, in section \ref{sec:bipartite} we give a completely new
application, to the ``scan'' chain for sampling colourings of
bipartite graphs.  For even relatively small values of $\Delta$, our
results improve Vigoda's~\cite{V99} celebrated $11\Delta/6$ bound on
the number of colours required for rapid mixing.

\section{Path coupling and stopping times}\label{sec:stopping}

We first deal with the most useful and applicable case, in which the stopping time for a pair of coupled chains is the first time that the distance between the two chains changes. This simplifies the proofs and makes the thrust of the argument clearer. In Section~\ref{sec:stopping2} we do deal with more general stopping times, however it should be noted that so far all applications of stopping times results in path coupling have only used this simple form of stopping time. 

\subsection{Distance-change stopping time}\label{sec:stopping1}
 
Let $\M$ be a Markov chain on state space $\Omega$. Let $\mathrm{d}$
be an integer valued metric on $\Omega\times \Omega$, and let
$(X_t,Y_t)$ be a path coupling for $\M$. We define $T_t$, a stopping
time for the pair $(X_t,Y_t)\in S$, to be the smallest $t'>t$ such that $\dd(X_{t'},Y_{t'})\neq \dd(X_t,Y_t)$. We will define a new metric
$\dd'$ such that if we have contraction in the metric $\dd$ at the
stopping times, then we have contraction in the metric $\dd'$ at
every step which has a positive probability of being a stopping
time. 

Let $\alpha>0$ be a constant such that $\E[\dd(X_{T_t},Y_{T_t})]\leq \alpha\dd(X_t,Y_t) \textrm{ for all
}(X_t,Y_t)\textrm{ in }S$. If $\alpha <1$, then for any $(X_t,Y_t)\in S$, we simply define $\dd'$ as follows.
\begin{equation}\label{medeq10}
\dd'(X_t,Y_t)= (1-\alpha)\dd(X_t,Y_t)+\E[\dd(X_{T_t},Y_{T_t})]\leq
\dd(X_t,Y_t).
\end{equation}
The metric is extended in the usual way to pairs $(X_t,Y_t)\notin
S$, using shortest paths. See, for example, \cite{DG99}. We will
apply path coupling with the metric $\dd'$ and the original
coupling. First we show a contraction property for this metric.
\begin{lem}\label{medlem1}
If $\E[\dd(X_{T_t},Y_{T_t})]\leq \alpha\dd(X_t,Y_t) <\dd(X_t,Y_t)\textrm{ for all
}(X_t,Y_t)\textrm{ in }S,$ then
\[
\E[\dd'(X_k,Y_k)\,|\, X_0,Y_0]\leq \big(1-(1-\alpha)\Pr(T_0\leq
k)\big)\dd'(X_0,Y_0).\]\end{lem}%
\begin{proof}%
We prove this by induction on $k$. It obviously holds for $k=0$,
since $T_0>0$. Using $\mathds{1}_{\mathcal{A}}$ to denote the 0/1
indicator of any event $\mathcal{A}$, we may write (\ref{medeq10})
as
\begin{align}
\dd'(X_0,Y_0)\ &=\ (1-\alpha)\dd(X_0,Y_0)+
\E[\dd(X_{T_k},Y_{T_k})\mathds{1}_{T_0>k}]+
\E[\dd(X_{T_0},Y_{T_0})\mathds{1}_{T_0\leq k}],\label{medeq20}
\end{align}
since if $T_0>k$ then $T_k=T_0$. Similarly, we have
\begin{align}
\E[\dd'(X_k,Y_k)]&\ =\ \E[\dd'(X_{k},Y_{k})\mathds{1}_{T_0>k}]+
\E[\dd'(X_{k},Y_{k})\mathds{1}_{T_0\leq k}]\notag\\
&\ =\ (1-\alpha)\E[\dd(X_{k},Y_{k})\mathds{1}_{T_0>k}]
+\E[\dd(X_{T_k},Y_{T_k})\mathds{1}_{T_0>k}]+
\E[\dd'(X_{k},Y_{k})\mathds{1}_{T_0\leq k}].\notag\\
&\ =\ (1-\alpha)\E[\dd(X_{0},Y_{0})\mathds{1}_{T_0>k}]
+\E[\dd(X_{T_k},Y_{T_k})\mathds{1}_{T_0>k}]+
\E[\dd'(X_{k},Y_{k})\mathds{1}_{T_0\leq k}].\label{medeq30}
\end{align}
Subtracting (\ref{medeq20}) from (\ref{medeq30}), we have
\begin{align*}
\E[\dd'(X_k,Y_k)]-\dd'(X_0,Y_0) &= \
-(1-\alpha)\E[\dd(X_{0},Y_{0})\mathds{1}_{T_0\leq k}]
+\E[(\dd'(X_k,Y_k)-\dd(X_{T_0},Y_{T_0}))\mathds{1}_{T_0\leq k}].
\end{align*}
For $T_0\leq k$, since $k-T_0\leq k-1$ the inductive hypothesis
implies $\E[\dd'(X_k,Y_k)\,|\, X_{T_0},Y_{T_0}] \leq
\dd'(X_{T_0},Y_{T_0})\leq \dd(X_{T_0},Y_{T_0})$, (if
$(X_k,Y_k)\not\in S$ this follows by linearity). Hence we have
\begin{align*}
\E[\dd'(X_k,Y_k)]-\dd'(X_0,Y_0) &\leq\
-(1-\alpha)\E[\dd(X_{0},Y_{0})\mathds{1}_{T_0\leq k}],
\end{align*}
The conclusion follows, since
$\E[\dd(X_{0},Y_{0})\mathds{1}_{T_0\leq k}]=\Pr(T_0\leq
k)\dd(X_{0},Y_{0})\geq \Pr(T_0\leq k)\dd'(X_{0},Y_{0})$.
\end{proof}

We may now prove the first version of our main result.
\begin{thm}\label{stopping1}
Let $\M$ be a Markov chain on state space $\Omega$. Let $\mathrm{d}$
be an integer valued metric on $\Omega$, and let $(X_t,Y_t)$ be a
path coupling for $\M$. Let $T_t$ be the above stopping times.
Suppose for all $(X_0,Y_0)\in S$ and for some integer $k$ and $p>0$,
that
\begin{enumerate}
\item[(i)] $\Pr[T_0 \leq k]\geq p$, %
\item[(ii)] $\E[\mathrm{d}(X_{T_0},Y_{T_0})/\dd(X_0,Y_0)]\leq \alpha  < 1$.
\end{enumerate}
Then the mixing time $\tau(\eps)$ of $\M$ satisfies
\[ \tau(\eps)\
\leq\frac{k}{p(1-\alpha)}\ln\Big(\frac{eD}{\eps(1-\alpha)}\Big).\]
where $D=\max\{\mathrm{d}(X,Y):X,Y\in\Omega\}$.
\end{thm}
\begin{proof}
>From Lemma \ref{medlem1}, $\dd'$ contracts by a factor
$1-(1-\alpha)p\leq e^{-(1-\alpha)p}$ for every $k$ steps of $\M$.
Note also that $\dd'\leq D$. It follows that, at time $\tau(\eps)$,
we have
\begin{equation*}
\Pr(X_\tau\neq Y_\tau)\leq \E[\dd(X_\tau,Y_\tau)]\leq
\frac{\E[\dd'(X_\tau,Y_\tau)]}{1-\alpha}\leq
\frac{De^{-(1-\alpha)p\tau/k}}{1-\alpha}\leq\eps,
\end{equation*}
from which the theorem follows.
\end{proof}
If $1-\alpha$ is small compared to $\eps$, it is possible to do
better than this. We will need the technical Lemma~\ref{medlm20}
below, which says that we will not have to wait too long for a
stopping time to occur. 


\begin{lem}\label{medlm20}
If $\M$ satisfies the conditions of Theorem~\ref{stopping1} 
then $\Pr[T_t>t+t']
\leq (1-p)^{\lfloor t'/k \rfloor}$.
\end{lem}
\begin{proof}
We prove this by induction on $t'$. It clearly holds for all $t$ and
$t'< k$ since $\lfloor t'/k \rfloor=0$. Suppose inductively that
$\Pr[T_t>s+t] \leq (1-p)^{\lfloor s/k \rfloor}$ for all $t$ and
$s<t'$. Then, if $t'\geq k$,
\begin{align*}
\Pr[T_t>t+t'] &=\ \Pr[T_{t}>t+t'-k \textrm{ and }
T_{t+t'-k}>t+t']\\
&=\ \Pr[T_{t}>t+t'-k] \Pr[T_{t+t'-k}>t+t'\ |\ T_t>t+t'-k].
\end{align*}
Since the process is Markovian, and by condition (i), 
\begin{align*}
\Pr[T_{t+t'-k}>t+t'\ |\ T_t>t+t'-k] &\leq\ \max\{\Pr[T_{t+t'-k}>t+t' ] :(X_{t+t'-k},Y_{t+t'-k})\in S\}\\
& =\ \max\{\Pr[T_{0}>k\ ]:(X_{0},Y_{0})\in S\}\\
&\leq\ 1-p.
\end{align*} By the inductive hypothesis this gives
\begin{align*}
\Pr[T_t>t'+t] &\leq \ (1-p)^{\lfloor (t'-k)/k \rfloor}(1-p)\ =\
(1-p)^{\lfloor t'/k \rfloor}.\qedhere
\end{align*}
\end{proof}
\begin{thm}\label{stopping}
Let $\M$ be a Markov chain on state space $\Omega$. Let $\mathrm{d}$
be an integer valued metric on $\Omega\times \Omega$, and let
$(X_t,Y_t)$ be a path coupling for $\M$. Let $T_t$ be the above stopping time.
Suppose for all $(X_0,Y_0)\in S$ and for some integer $k$ and $p>0$,
that
\begin{enumerate}
\item[(i)] $\Pr[T_0\leq k ]\geq p,$
\item[(ii)] $\E[\mathrm{d}(X_{T_0},Y_{T_0})/\dd(X_0,Y_0)]\leq \alpha  < 1$.
\end{enumerate}
Then the mixing time
$\tau(\eps)$ of $\M$ satisfies
\[ \tau(\eps)\leq
\frac{k(2-\alpha)}{p(1-\alpha)}\ln\Big(\frac{2eD}{\eps}\Big).
\]
where $D=\max\{\mathrm{d}(X,Y):X,Y\in\Omega\}$.
\end{thm}
\begin{proof}
Let $X_t=Z^0_0, Z^1_0,\ldots,Z^r_0=Y_t$ be a shortest path from
$X_t$ to $Y_t$ in the metric $\dd'$, such that $(Z^i_0,Z^{i+1}_0)\in
S$ ($i=0,\ldots,r-1$). If $\mathfrak{t}_i$ is the stopping  time for
$(Z^i_0,Z^{i+1}_0)$ then, using Lemma \ref{medlm20},
\begin{align*}
\Pr(X_{t+t'}\neq Y_{t+t'}\,|\,X_{t},Y_{t})&\leq\ \Pr(\exists i : Z^i_{t'}\neq Z^{i+1}_{t'})\\
&\leq\ \Pr(\exists i : Z^i_{\mathfrak{t}_i}\neq
Z^{i+1}_{\mathfrak{t}_i}\textrm{ or }\mathfrak{t}_i>t')\\
&\leq\ \sum_{i=0}^{r-1} \left( \E[\dd(Z^i_{\mathfrak{t}_i},
Z^{i+1}_{\mathfrak{t}_i})]
+ \Pr(\mathfrak{t}_i>t')\right)\\
&\leq\ \sum_{i=0}^{r-1} \left( \dd'(Z^i_0, Z^{i+1}_0) + (1-p)^{\lfloor t'/k \rfloor} \right)\\
&\leq\  \dd'(X_t, Y_t) + D(1-p)^{\lfloor t'/k \rfloor}.
\intertext{Hence}\ \Pr(X_{t+t'}\neq Y_{t+t'})&\leq \E[\dd'(X_t,
Y_t)]
+ D(1-p)^{\lfloor t'/k \rfloor},\\
&\leq\ De^{-(1-\alpha)p\lfloor t/k\rfloor} + D(1-p)^{\lfloor t'/k
\rfloor}.\\
&\leq\ D\big(e^{-(1-\alpha)p\lfloor t/k\rfloor} + e^{-p\lfloor
t'/k\rfloor}\big).\\
\intertext{Therefore} \Pr(X_{t+t'}\neq Y_{t+t'})\ &\leq\
\tfrac{1}{2}\eps+\tfrac{1}{2}\eps\ =\ \eps,\ \ \mathrm{if} \ t\geq
k\Big\lceil\frac{\ln(2D/\eps)}{p(1-\alpha)}\Big\rceil\ \mathrm{and}
\ t'\geq k\Big\lceil\frac{\ln(2D/\eps)}{p}\Big\rceil.
\end{align*}
The statement of the theorem now follows easily.
\end{proof}

\subsection{General stopping times}\label{sec:stopping2}

We now extend the results proved in this section to incorporate stopping times other than the first time at which the distance changes. In order to make sense in the context of path coupling, the stopping times must satisfy the following conditions.

\begin{minipage}[l]{0.9\textwidth}\vspace{2ex}
\textsc{Stopping time conditions:}
\begin{enumerate}
\item There must be a stopping time $T(X_0,Y_0)$ defined for each pair $(X_0,Y_0)\in S$, such\\ that $\E[\dd(X_{T_0(X,Y)}, Y_{T_0(X,Y)})]\leq \alpha \dd(X_0,Y_0)$. \label{contraction}
\item For all $(X_0,Y_0)\in S$ we have $\Pr[T(X_0,Y_0)\leq k ]\geq p$.\label{shorttime}
\item The coupling should be Markovian. \label{markovian}
\end{enumerate}
\end{minipage}\vspace{1ex}

We may assume that for $(X_0, Y_0)\in S$ if $X_t=Y_t$ then $T(X_0,Y_0)\leq t$. Since the future evolution of $(X_t,Y_t)$ does not depend on the evolution up to time $t$, by \ref{contraction} and \ref{markovian} it follows that for all $t\geq 0$ there is a stopping time $T_t(X,Y)$ such that if $X_t=X, Y_t=Y$ then $\E[\dd(X_{T_t(X,Y)}, Y_{T_t(X,Y)})]\leq \alpha \dd(X_t,Y_t)$. Moreover, from \ref{shorttime} and \ref{markovian}, it follows that $\Pr[T_t(X,Y)\leq k+t ]\geq p$. 

When dealing with the first change in distance we had the benefit that for all $(X_t,Y_t)\in S$ and $t'>t$, if $T_t(X_t,Y_t)>t'$ then $(X_{t'},Y_{t'})\in S$  and also $T_{t'}(X_{t'},Y_{t'})=T_t(X_t,Y_t)$. This no longer necessarily holds. We must therefore be more careful about exactly which stopping time we are referring to at any time and regarding any pair of states.

Let $(X_t,Y_t)$ be a coupled evolution of the chain, and let $P_t=(X_t=Z^0_t, Z^1_t,Z^2_t,\ldots,Z^{d_t}_t=Y_t)$ be the path-coupling path from $X_t$ to $Y_t$, so that $(Z^i_t,Z^{i+1}_t)\in S$ for all $i,t$. We will inductively define a set of \emph{starting pairs} in the paths $P_t$, $t\geq0$ as follows. 
\begin{enumerate}
\item For all $i, (Z^i_0,Z^{i+1}_0)$ is a starting pair. 
\item For each $(Z^i_{t_1},Z^{i+1}_{t_1})$ if there is a time $t_0\leq t_1$ and {starting pair} $(Z^j_{t_0}, Z^{j+1}_{t_0})\in P_{t_0}$ such that $(Z^i_{t_1},Z^{i+1}_{t_1})$ is in the subpath of $P_{t_1}$ which evolved from $(Z^j_{t_0}, Z^{j+1}_{t_0})$ and $T_{t_0}(Z^j_{t_0}, Z^{j+1}_{t_0})>t_1$, then $(Z^i_{t_1},Z^{i+1}_{t_1})$ is \emph{not} a starting pair but $t_0$ is the \emph{starting time associated with} $(Z^i_{t_1},Z^{i+1}_{t_1})$ and $(Z^j_{t_0}, Z^{j+1}_{t_0})$ is the \emph{starting pair associated with }$(Z^i_{t_1},Z^{i+1}_{t_1})$. 
\item For each $(Z^i_{t_1},Z^{i+1}_{t_1})$ such that there is no time and pair as above, then $(Z^i_{t_1},Z^{i+1}_{t_1})$ is defined to be a starting pair. Note that in this case there must be a time $t_0\leq t_1$ and {starting pair} $(Z^j_{t_0}, Z^{j+1}_{t_0})\in P_{t_0}$ such that $(Z^i_{t_1},Z^{i+1}_{t_1})$ is in the subpath of $P_{t_1}$ which evolved from $(Z^j_{t_0}, Z^{j+1}_{t_0})$ and $T_{t_0}(Z^j_{t_0}, Z^{j+1}_{t_0})=t_1$. 
\end{enumerate}
For a starting pair $(Z^j_{t_0}, Z^{j+1}_{t_0})$, we define the distance at time $t_1, t_0\leq t_1<T_{t_0}(Z^j_{t_0}, Z^{j+1}_{t_0})$ to be 
\begin{align} \label{startdist}
d_{t_1}(Z^j_{t_0}, Z^{j+1}_{t_0})= (1-\alpha)\dd(Z^j_{t_0}, Z^{j+1}_{t_0}) +\E\left[\dd(Z^j_{T_{t_0}(Z^j_{t_0} Z^{j+1}_{t_0})}, Z^{j+1}_{T_{t_0}(Z^j_{t_0} Z^{j+1}_{t_0})})\ |\ \mathcal{F}_{t_1}\right]
\end{align}
where $\mathcal{F}_t$ is the
$\sigma$-algebra generated by $\{(X_{t'},Y_{t'}):t'\leq t\}$. Thus
$\{\mathcal{F}_t:t'\geq 0\}$ is the filtration generated by the
coupling. The distance at times not in the given range is zero. This is analogous to the definition of the new metric in equation~(\ref{medeq10}). At a time $t$ we are interested in the set $\mathcal{SP}_t$ of starting pairs $(Z^j_{t_0}, Z^{j+1}_{t_0})$ for which $T_{t_0}(Z^j_{t_0}, Z^{j+1}_{t_0})> t$. We define the distance between $X_t$ and $Y_t$ to be 
\begin{align}
D(X_t,Y_t)=\sum_{(Z^j_{t_0}, Z^{j+1}_{t_0})\in \mathcal{SP}_t} d_{t}(Z^j_{t_0}, Z^{j+1}_{t_0}). \label{dt}
\end{align}
It is clear that if $\dd(X_t,Y_t)\neq 0$ then $D(X_t,Y_t)\geq (1-\alpha)$. 
We now prove a contraction lemma analagous to Lemma~\ref{medlem1}.
\begin{lem}
Given the stopping times conditions, then for all $(X_0,Y_0)$ and all $t\geq 0$
\[
\E[D(X_{t+k},Y_{t+k})\,|\, X_t,Y_t]\leq \left(1-\frac{(1-\alpha)p}{\gamma+1}\right)D(X_k,Y_k),\]
where $\gamma$ is the maximum value of $\E[\dd(X_{T_0(X,Y)}, Y_{T_0(X,Y)})\ |\ \mathcal{F}_t ]/\dd(X_0,Y_0)$ over all pairs in $S$ and evolutions $\mathcal{F}_t$ such that $t< T_0(X,Y)$. 
\end{lem}%
\begin{proof}%
The set $\mathcal{SP}_{t+k}$ is the union of the starting pairs from $\mathcal{SP}_{t}$ which did not reach their stopping time by time $t+k$, and those starting pairs arising from a pair in $\mathcal{SP}_{t}$ which did stop by time $t+k$. Hence, writing $T_{t_0}$ for ${T_{t_0}(Z^j_{t_0} Z^{j+1}_{t_0})}$, 
\begin{align*}
D(X_{t+k},Y_{t+k})\ &=\ \sum_{(Z^j_{t_0}, Z^{j+1}_{t_0})\in \mathcal{SP}_t} \mathds{1}_{T_{t_0}> t+k}d_{t+k}(Z^j_{t_0}, Z^{j+1}_{t_0})  + \mathds{1}_{T_{t_0}\leq t+k}\sum d_{t+k}(Z^l_{t_l}, Z^{l+1}_{t_l})
\end{align*}
where the second sum is over starting pairs arising from the stopping of pair $(Z^j_{t_0}, Z^{j+1}_{t_0})$. As in Lemma~\ref{medlem1}, we may assume inductively that $\E[d_{t+k}(Z^l_{t_l}, Z^{l+1}_{t_l})\ | \ t< t_l\leq t+k]\leq \dd(Z^l_{t_l}, Z^{l+1}_{t_l})$. Then, given $\mathcal{F}_t$, the expected value of $D(X_{t+k},Y_{t+k})$ is
\begin{align}
\E[D(X_{t+k},Y_{t+k})]\ &\leq\ \sum_{(Z^j_{t_0}, Z^{j+1}_{t_0})\in \mathcal{SP}_t} \E[\mathds{1}_{T_{t_0}> t+k}d_{t+k}(Z^j_{t_0}, Z^{j+1}_{t_0})]  + \E[\mathds{1}_{T_{t_0}\leq t+k}\dd(Z^j_{T_{t_0}}, Z^{j+1}_{T_{t_0}})]\notag \\
 &\leq\ \sum_{(Z^j_{t_0}, Z^{j+1}_{t_0})\in \mathcal{SP}_t} \E\left[\mathds{1}_{T_{t_0}>t+k} (1-\alpha)\dd(Z^j_{t_0}, Z^{j+1}_{t_0})\right] +\E[\dd(Z^j_{T_{t_0}}, Z^{j+1}_{T_{t_0}})].\label{tplusone}
\end{align}
So subtracting (\ref{dt}) from (\ref{tplusone}) we get
\begin{align}
\E[D(X_{t+k},Y_{t+k})]-D(X_t,Y_t) &\leq \ \sum_{(Z^j_{t_0}, Z^{j+1}_{t_0})\in \mathcal{SP}_t}
-(1-\alpha)\E[\dd(Z^j_{t_0}, Z^{j+1}_{t_0})\mathds{1}_{T_0\leq t+k}]\notag \\
& \leq \   -(1-\alpha)p \sum_{(Z^j_{t_0}, Z^{j+1}_{t_0})\in \mathcal{SP}_t} \dd(Z^j_{t_0}, Z^{j+1}_{t_0})\label{wrongcontraction} \\
&\leq \ -\frac{(1-\alpha)p}{1-\alpha+\gamma}  D(X_t,Y_t). \notag
\end{align}
The final inequality follows since, by (\ref{startdist}), we have $ d_{t}(Z^j_{t_0}, Z^{j+1}_{t_0})\leq (1-\alpha+\gamma) \dd(Z^j_{t_0}, Z^{j+1}_{t_0}).$
\end{proof}
The $\gamma$ term arises because although we have contraction in inequality~(\ref{wrongcontraction}), we need to express this as a proportion of $D(X_t,Y_t)$. The expected value at the stopping time is only guaranteed to be at most $\alpha \dd$ at the outset. If we have already evolved, possibly adversely, the expected value at the stopping time could be larger than this, and the proportional changes correspondingly smaller. However $\gamma$ is bounded by the maximum distance (in the original metric) that can occur at the stopping time; in practice this is very likely to be a small constant.

By following the same arguments as in Section~\ref{sec:stopping1}, with this contraction lemma we obtain the following theorem.
\begin{thm}\label{stoppinggeneral}
Let $\M$ be a Markov chain on state space $\Omega$. Let $\mathrm{d}$
be an integer valued metric on $\Omega$, and let $(X_t,Y_t)$ be a
path coupling for $\M$. Let $T(X_0,Y_0)$ be stopping times satisfying the stopping times conditions.
Then the mixing time $\tau(\eps)$ of $\M$ satisfies

\[ \tau(\eps)= O\left(\frac{k(1-\alpha+\gamma)}{p(1-\alpha)}\ln\Big(\frac{D}{\eps}\Big)\right).
\]
\end{thm}

\begin{rem}
One of the most interesting features of Theorems~\ref{stopping} and~\ref{stoppinggeneral} is that their
proofs employ only standard path coupling (applied to the $k$-step chain), but with a metric  which
has some useful properties. Thus, for any problem to which stopping
times might be applied, there exists a metric from which the same
result could be obtained using one-step path coupling.
\end{rem}

\begin{rem}\label{rem15}
Stopping times condition~\ref{shorttime} may appear a restriction, but appears to be
naturally satisfied in most applications, even with $k=1$. The
alternative, though less natural, assumption of uniformly bounded
stopping times~\cite{HV04} is also included. (See Remark
\ref{rem20}.)
\end{rem}

\begin{rem}\label{rem20}
We may compare this stopping time theorem with those
in~\cite{BDK05a,HV04}. The main result of~\cite[Theorem 3]{HV04}
concerns bounded stopping times, where $T_0\leq M$ for all
$(X_0,Y_0)\in S$, and gives a mixing time of $O(M
(1-\alpha)^{-1}\log D).$ By setting $k=M$ and $p=1$ in
Theorem~\ref{stopping}, we obtain the same mixing time up to minor
changes in constants, but with a proof that does not involve
defining a multistep coupling. For unbounded mixing times,
\cite[Corollary~4]{HV04} gives a bound $O(\E[T] (1-\alpha)^{-2}W\log
D)$ by truncating the stopping times, where $W$ denotes the maximum
of $\dd(X_t,Y_t)$ over all $(X_0, Y_0) \in S$ and $t\leq T$. In most
applications $\E[T]\leq k/p$,  so in Theorem~\ref{stopping} we obtain an improvement of order
$W(1-\alpha)^{-1}$. By comparison with~\cite{BDK05a},  we obtain a
more modest improvement, of order $\log W\log(D(1-\alpha)^{-1})/\log
D$. For the more general stopping times, comparing Theorem~\ref{stoppinggeneral} and \cite[Corollary~4]{HV04}, we obtain an improvement of order $\frac{W}{\gamma(1-\alpha)}$. It should be noted that $\gamma\leq W$. 
\end{rem}
\begin{rem}\label{rem25}
Further improvements to Theorem~\ref{stopping} seem unlikely, other
than in constants. The term  $k/p$ must be present, since it bounds
a single stopping time. A term $1/(1-\alpha)\log(D/\varepsilon)=
\Theta(\log_\alpha(D/\varepsilon))$ also seems essential, since it
bounds the number of stopping times required. Likewise improvements to Theorem~\ref{stoppinggeneral} are likely restricted to changing the dependence on $\gamma$, although it seems plausible that some dependence is required.
\end{rem}

\section{Hypergraph independent sets}\label{sec:indsets}
We now turn our attention to hypergraph independent sets. These were
previously studied in~\cite{BDK05a}. Let $\CH=(\CV,\CE)$ be a
hypergraph of maximum degree $\Delta$ and minimum edge size $m$. A
subset $S\subseteq \CV$ of the vertices is \emph{independent} if no
edge is a subset of $S$. Let $\Omega(\CH)$ be the set of all
independent sets of $\CH$. We define the Markov chain $\M(\CH)$ with
state space $\Omega(\CH)$ by the following transition process
(\emph{Glauber dynamics}). If the state of $\M$ at time $t$ is
$X_t$, the state at $t+1$ is determined by the following procedure.
\begin{enumerate}
\item Select a vertex $v\in \CV$ uniformly at random, \item
\begin{enumerate}
    \item if $v\in X_t$ let $X_{t+1}=X_t\backslash \{v\}$ with probability
    $1/2$,
    \item if $v\not\in X_t$ and $X_t \cup \{v\} $ is independent, let
    $X_{t+1}=X_t\cup \{v\}$ with probability $1/2$,
    \item otherwise let $X_{t+1}=X_t$.
\end{enumerate}
\end{enumerate}
This chain is easily shown to be ergodic with uniform stationary
distribution. The natural coupling for this chain is the
``identity'' coupling, the same transition is attempted in both
copies of the chain. If we try to apply standard path coupling to
this chain, we immediately run into difficulties. Consider a state
of the coupled chain at a time $t$, $(X_t,Y_t)$, such that
$Y_t=X_t\cup\{w\}$, where $w\notin X_t$ (the \emph{change vertex})
is of degree $\Delta$. An edge $e\in\CE$ is \emph{critical} in $Y_t$
if it has only one vertex $z\in\CV$ which is not in $Y_t$, and we
call $z$ \emph{critical for $e$}. If each of the edges through $w$
is critical for $Y_t$, then there are $\Delta$ choices of $v$ in the
transition which can be added in $X_t$ but not in $Y_t$. Thus the
change in the expected \emph{Hamming distance} between $X_t$ and
$Y_t$ after one step could be as high as
$\frac{\Delta}{2n}-\frac1n$, and we obtain rapid mixing only in the
case $\Delta=2$.

For $(\sigma, \sigma\cup \{w\})\in S$, let $E_i(w,\sigma)$ be the
set of edges containing $w$ which have $i$ occupied vertices in
$\sigma$. Using a result like Theorem~\ref{stopping1} above, it is
shown in~\cite{BDK05b} that, for the stopping time $T$ given by the
first epoch at which the Hamming distance between the coupled chains
changes,
\begin{equation*}
\E[\dd_{\mathrm{Ham}}(X_T,Y_T|X_0=\sigma, Y_0=\sigma\cup \{w\})]\leq
2\sum_{i=0}^{m-2} p_i|E_i|\leq 2p_1\Delta,
\end{equation*}
where the $p_i$ is the probability that $\dd(X_T,Y_T)=2$ if $w$ is
in a single edge with $i$ occupied vertices. Since  $p_1< 1/(m-1)$,
we  obtain rapid mixing  when  $2\Delta/(m-1)\leq 1$, i.e. when
$m\geq 2\Delta+1$. See~\cite{BDK05b} for details.

The approach of Section~\ref{sec:stopping} would lead us to define a
metric for which the distance between $\sigma$ and $\sigma\cup
\{w\}$ is $(1-2 p_1\Delta) + 2\sum_{i=0}^{m-2} p_i|E_i|.$ By
Lemma~\ref{medlem1}, we know that this metric contracts in
expectation. However, prompted by the form of this metric, but
retaining the freedom to optimise constants, we will instead define
the new metric $\dd$ to be
\[
\dd(\sigma, \sigma\cup \{w\})=\sum_{i=0}^{m-2} c_i|E_i| ,
\]
where $0<c_i\leq 1$ $(0\leq i\leq m-2)$ are a nondecreasing sequence
of constants to be determined. Using this metric, we obtain the
following theorem.
\begin{thm}\label{indsets}
Let $\Delta$ be fixed, and let $\CH$ be a hypergraph such that
$m\geq \Delta+2\geq 5$, or $\Delta=3$ and $m\geq 2$. Then the Markov
chain $\M(\CH)$ has mixing time $O(n\log n)$.
\end{thm}

\begin{proof}
Without loss of generality, we take $c_{m-2}=1$ and we will define
$c_{-1}=c_0, c_{m-1}\geq\Delta+1$. Note that $c_{-1}$ has no real
role in the analysis, and is chosen only for convenience, but
$c_{m-1}$ is chosen so that $c_{m-1}-c_{m-2} \geq\Delta \geq
\dd(\sigma,\sigma')$ for any pair $(\sigma,\sigma')\in S$. We
require $c_i>0$ for all $i$ so that we will always have $\dd(\sigma,
\sigma')>0$ if $\sigma\neq \sigma'$.

Now consider the expected change in distance between $\sigma$ and
$\sigma\cup \{w\}$ after one step of the chain.

If $w$ is chosen, then the distance decreases by $\sum_{i=0}^{m-2}
c_i|E_i| $. The contribution to the expected change in distance is
$- \frac{2}{2n}\sum_{i=0}^{m-2}  c_i|E_i|$.

If we insert a vertex $v$ in an edge containing $w$, then we
increase the distance by $(c_{i+1}-c_i)\geq 0$ for each edge in
$E_i$ containing $v$. This holds for $i=0,\ldots,m-2$, by the choice
of $c_{m-1}=\Delta+1$. Let $U$ be the set of unoccupied neighbours
of $w$, and $\nu_i(v)$ be the number of edges with $i$ occupants
containing $w$ and $v$. Then the contribution is
\[
 \sum_{v\in U} \frac 1{2n} \sum_{i=0}^{m-2} \nu_i(v) (c_{i+1}-c_i) =
 \frac 1{2n} \sum_{i=0}^{m-2}  (c_{i+1}-c_i)(m-i-1)|E_i|,
\]
since
\begin{equation*}
    \sum_{v\in U} \nu_i(v) = \sum_{v\in U} \sum_{e\in E_i} \mathds{1}_{v\in e}
    = \sum_{e\in E_i} \sum_{v\in e\cap U} 1 = \sum_{e\in E_i} (m-i-1) = (m-i-1)|E_i|.
\end{equation*}

If we delete a vertex $v$ in an edge containing $w$, then we
decrease the distance by $(c_{i}-c_{i-1})$ for each edge in $E_i$
containing $v$. This holds for $i=0,\ldots,m-2$, by the choice of
$c_{-1}$. Let $O$ be the set of occupied neighbours of $w$, and
$\nu_i(v)$ be the number of edges with $i$ occupants containing $w$
and $v$. Then the contribution is
\[
 - \sum_{v\in O} \frac 1{2n} \sum_{i=0}^{m-2} \nu_i(v) (c_{i}-c_{i-1})
 = - \frac 1{2n} \sum_{i=0}^{m-2} (c_{i}-c_{i-1})i|E_i|,
\]
since, as for $U$ above,
\begin{equation*}
    \sum_{v\in O} \nu_i(v) = \sum_{v\in O} \sum_{e\in E_i} \mathds{1}_{v\in e}
    = \sum_{e\in E_i} \sum_{v\in e\cap O} 1 = \sum_{e\in E_i} i = i|E_i|.
\end{equation*}

Let $\dd_0=\dd(\sigma, \sigma\cup \{w\})$, and let $\dd_1$ be the
distance between the evolved states after one step of the chain. The
change in expected distance $\E[\dd_1-\dd_0]$ satisfies

\begin{align*}
2n \E[\dd_1-\dd_0] &\leq\  - 2 \sum_{i=0}^{m-2} c_i|E_i|  +
\sum_{i=0}^{m-2} (c_{i+1}-c_i)(m-i-1)|E_i|- \sum_{i=0}^{m-2} (c_{i}-c_{i-1})i|E_i| \\
&=\  \sum_{i=0}^{m-2}  \left( -2 c_i +(m-i-1)(c_{i+1}-c_i)- i(c_{i}-c_{i-1})\right)|E_i|\\
&=\  \sum_{i=0}^{m-2}  \left( ic_{i-1}-(m+1)c_i+(m-i-1)c_{i+1}
\right)|E_i|.
\end{align*}

We require $ \E[\dd_1-\dd_0]\leq -\gamma$, for some $\gamma\geq 0$,
which holds for all possible choices of $E_i$ if and only if
$(m-i-1)c_{i+1}-(m+1)c_i +ic_{i-1}\leq -\gamma$ for all
$i=0,1,\ldots,m-2$.  Thus we need a solution to
\begin{align}
&ic_{i-1}-(m+1)c_i+(m-i-1)c_{i+1}\leq -\gamma\qquad(i=0,\ldots,m-2),\label{eq10}\\
&0=c_{-1}< c_{0}\leq c_{1}\leq \cdots \leq c_{m-3}\leq c_{m-2}= 1, \notag\\
&c_{m-1}\geq \Delta+1,\ \gamma\geq 0,\notag
\end{align}
with $\gamma>0$ if possible. Adding~(\ref{eq10}) from $i$ to $m-2$
gives
\begin{align}
  &ic_{i-1}-(m-i)c_i-(m-1)c_{m-2}+c_{m-1}\ \leq\ -(m-i-1)\gamma&(i=0,\ldots,m-2),\notag\\
\mathrm{i.e.}\ \ &ic_{i-1}\ \leq\
(m-i)c_i+(m-\Delta-2)-(m-i-1)\gamma &(i=0,\ldots,m-1).\label{eq20}
\end{align}

Substitute $u_i=\binom{m-1}{i}c_i$ in (\ref{eq20}), so
$u_{m-1}\geq\Delta+1$, $u_{m-2}=m-1$ and $u_{-1}=0$. Then we have
\begin{equation*}
u_{i-1}\leq u_i+\frac{m-\Delta-2+\gamma}{m}\binom{m}{i}
-\gamma\binom{m-1}{i}\quad \quad(i=0,\ldots,m-2).
\end{equation*}
Using the boundary condition $u_{-1}=0$, these give
\begin{equation*}
u_i \leq \gamma\sum_{j=0}^{i}\binom{m-1}{i}
-\frac{m-\Delta-2+\gamma}{m}\sum_{j=0}^{i}\binom{m}{j}
\quad(i=0,\ldots,m-2).
\end{equation*}
The boundary condition $u_{m-2}=m-1$ now implies
\[\gamma\ \leq\ \frac{2^{m}-1-m}{(m-2)2^{m-1}+1}
\left(m-\Delta-2+\frac{m(m-1)}{2^{m}-1-m}\right).\] Let
\[ f(m) = m-2 +\frac{m(m-1)}{2^m-1-m},\]
then we can have $\gamma\geq 0$ if and only if $f(m)\geq\Delta$, and
$\gamma>0$ if and only if $f(m)>\Delta$. Then
 \[c_i\ = \ \frac{\gamma\sum_{j=0}^{i}\binom{m-1}{j}
-\frac{m-\Delta-2+\gamma}{m}\sum_{j=0}^{i}\binom{m}{j}}
{\binom{m-1}{i}}\qquad(i=0,\ldots,m-2).\] In order to satisfy the
conditions of~(\ref{eq10}), we need to establish that $0< c_i\leq
c_{i+1}$ $(i=0,\ldots,m-3)$.

{\allowdisplaybreaks \begin{align*}
   c_i\ &= \ \frac{\gamma\sum_{j=0}^{i}\binom{m-1}{j}
-\frac{m-\Delta-2+\gamma}{m}\sum_{j=0}^{i}\binom{m}{j}}
{\binom{m-1}{i}}\qquad(i=0,\ldots,m-2)\\
&= \ \gamma\frac{\sum_{j=0}^{i}\binom{m-1}{j}
-\kappa\sum_{j=0}^{i}\binom{m}{j}} {\binom{m-1}{i}},\
\mathrm{where}\ \ \kappa\ =\ \frac{m-\Delta-2+\gamma}{m\gamma}.\\
&= \ \gamma\frac{\sum_{j=0}^{i}\binom{m-1}{j}
-\kappa\sum_{j=0}^{i}\left(\binom{m-1}{j}+\binom{m-1}{j-1}\right)} {\binom{m-1}{i}},\\
&= \ \gamma\frac{\sum_{j=0}^{i}\binom{m-1}{j}
-\kappa\left(\sum_{j=0}^{i}\binom{m-1}{j}+\sum_{j=0}^{i-1}\binom{m-1}{j}\right)} {\binom{m-1}{i}},\\
&= \ \gamma\frac{\sum_{j=0}^{i}\binom{m-1}{j}
-\kappa\left(2\sum_{j=0}^{i}\binom{m-1}{j}-\binom{m-1}{i}\right)} {\binom{m-1}{i}},\\
&= \ \gamma\frac{(1-2\kappa)\sum_{j=0}^{i}\binom{m-1}{j}
+\kappa\binom{m-1}{i}} {\binom{m-1}{i}},\\
&= \
\gamma(1-2\kappa)\frac{\sum_{j=0}^{i}\binom{m-1}{j}}{\binom{m-1}{i}}+\gamma\kappa,\\
&= \ \gamma(1-2\kappa) g_{i}+\gamma \kappa,\ \ \mathrm{say}.
\end{align*}}
Now $2\kappa< 1$ is equivalent to $2(m-\Delta-2)/(m-2)< \gamma$,
i.e.
\[ \frac{2(m-\Delta-2)}{m-2} < \frac{(2^m-1-m)(m-\Delta-2)+m(m-1)}{(m-2)2^{m-1}+1},\]
which holds for all $\Delta> 0$. Finally, $g_{i}$ is strictly
increasing, since
\begin{align*}
    \frac{g_{i-1}}{g_{i}} \ &=\
    \frac{\frac{m-i}{i}\sum_{j=0}^{i-1}\binom{m-1}{j}}{\sum_{j=0}^{i}\binom{m-1}{j}}\\
    \ &=\
    \frac{\sum_{j=1}^{i}\frac{m-i}{i}\binom{m-1}{j-1}}{\sum_{j=0}^{i}\binom{m-1}{j}}\\
    \ &\leq\
    \frac{\sum_{j=1}^{i}\binom{m-1}{j}}{\sum_{j=0}^{i}\binom{m-1}{j}},\ \ \mathrm{since}\
    j\leq i,\\
    \ &<\ 1.
\end{align*}
Hence $c_i$ is strictly increasing. It  only remains to verily that
$c_0>0$. This is clearly equivalent to $\gamma >
(m-\Delta-2)/(m-1)$. If $m=\Delta+2$, it follows from $\gamma> 0$.
If $m>\Delta+2$, it follows from $\gamma> 2(m-\Delta-2)/(m-2)$,
which we have already established.

If $m\geq 5$ then $m(m-1)/(2^m-1-m)<1$, so we will have $f(m)>
\Delta$ exactly when $m\geq\Delta+2$. For smaller values of $m$,
\begin{equation*}\renewcommand{\arraystretch}{1.5}
    \begin{array}{|c||c|c|c|}\hline
       m & 2 & 3 & 4\\\hline
       f(m) &\ 2\ & \,2\tfrac{1}{2}\, &
       3\tfrac{1}{11}\\\hline
     \end{array}
\end{equation*}
The new case here is $\Delta=3,m \geq 4$. In any case for which
$f(m)> \Delta$, standard path coupling arguments yield the mixing
times claimed since we have contraction in the metric and the
minimum distance is at least $c_0$. Since we can show mixing for
$\Delta=3,m\leq 3$ by other means (see~\cite{DG00a}), we have mixing
for $\Delta=3$ and every $m$.\end{proof}
\begin{rem}
The independent set problem here has a natural \emph{dual}, that of
sampling an \emph{edge cover} from a hypergraph with edge size
$\Delta$ and degree $m$. An edge cover is a subset of $\mathcal{E}$
whose union contains $V$. For the graph case of this sampling
problem, with arbitrary $m$, see~\cite{BD97a}. By duality this gives
the case $\Delta=2$ of the independent set problem here.
\end{rem}

\section{Satisfying assignments of SAT instances}\label{sec:sat}

The set of {\em independent} sets in a hypergraph with edge size $m$
and degree $\Delta$ corresponds in a natural way to the set of
{\em satisfying} assignments in a SAT instance with clause size equal to $m$
and number of each variable occurrences bounded by $\Delta$, cf. \cite{DG00a}.
The optimisation problems connected to small (variable) occurrence number instances of SAT were
studied recently in \cite{BKS03} (see \cite{BKS03} also for
additional references).

Given a hypergraph $\CH=(\CV,\CE)$ with $n$ vertices, $k$ hyperedges,
and edge size $m$ and degree $\Delta$. We construct an $m$SAT  formula
$f$, over $n$ variables $X=\{x_1,\ldots,x_n\}$ corresponding to vertices of
$\CH$ as follows.
If $e=\{v_1,\ldots,v_m\}$ is an hyperedge of $\CH$, we associate with
$e$ a clause $C_e=\bigvee^m_{i=1}\bar{x}_i$, and furthermore we set
$f=\bigwedge_{e\in\CE}C_e$. Notice that the number of satisfying assignments
of $f$ is precisely the same as a number of all independent sets of $\CH$,
and a number of occurrences of variables in $f$ is less than or equal to the degree
of $\CH$. We can moreover replace the literals $\bar{x}_i$ by $x_i$,
to obtain a monotone $m$SAT formula $f'$ with the same number of variable
occurrences as $f$ and with the same number of satisfying assignments.
The above construction is reversable, showing the
equivalence of corresponding counting problems of hypergraph independent
sets and monotone SAT formulas.

Let us denote by $\#(m,\Delta)\mu$SAT the problem of counting number of
satisfying assignments in monotone $m$SAT instances with at most $\Delta$
variable occurrences. Theorem \ref{indsets}. yields the first FPRASs
({\em Fully Polynomial Randomized Approximation Schemes}) for a large class
of monotone $m$SAT formulas.

\begin{thm}
  Let $\Delta$ be fixed, and $m\geq\Delta+2\geq 5$, or if $\Delta=3$ then
  $m\geq 2$. Then the associated Markov chain $\M(\CH)$ yields an FPRAS
  for the $\#(m,\Delta)\mu$SAT problem.
\end{thm}

The above result improves vastly the hitherto known results for approximate
counting the number of satisfying assignments of general monotone SAT formulas.

\section{Colouring 3-uniform hypergraphs}\label{sec:3uniform}

In our second application, also from~\cite{BDK05a}, we consider
proper colourings of 3-uniform hypergraphs. We again use Glauber
dynamics. Our hypergraph $\CH$ will have maximum degree $\Delta$,
uniform edge size $3$, and we will have a set of $q$ colours. For a
discussion of the easier problem of colouring hypergraphs with
larger edge size see~\cite{BDK05b}. A colouring of the vertices of
$\CH$ is proper if no edge is monochromatic. Let $\Omega'(\CH)$ be
the set of all proper $q$-colourings of $\CH$. We define the Markov
chain $\C(\CH)$ with state space $\Omega'(\CH)$ by the following
transition process. If the state of $\C$ at time $t$ is $X_t$, the
state at $t+1$ is determined by
\begin{enumerate}
\item selecting a vertex $v\in \CV$ and a colour $k\in\{1,2,\ldots,q\}$
uniformly at random, \item let $X'_t$ be the colouring obtained by
recolouring $v$
colour $k$ \item if $X'_t$ is a proper colouring let $X_{t+1}=X'_t$\\
otherwise let $X_{t+1}=X_t$.
\end{enumerate}
This chain is easily shown to be ergodic with the uniform stationary
distribution. For some large enough constant $\Delta_0$, it was
shown in~\cite{BDK05b} to be rapidly mixing for $q>1.65\Delta$ and
$\Delta>\Delta_0$, using a stopping times analysis. Here we improve
this result, and simplify the proof, by using a carefully chosen
metric which is prompted by the new insight into stopping times
analyses. If $w$  is the change vertex, the intuition
in~\cite{BDK05b} was that edges which contain both colours of $w$
are initially ``dangerous'' but tend to become less so after a time.
Thus our metric will be a function of the numbers of edges
containing $w$ with various relevant colourings.

\begin{thm}\label{cols}
Let $\Delta$ be fixed, and let $\CH$ be a 3-uniform hypergraph of
maximum degree $\Delta$. Then if $q \geq \big\lceil
\tfrac{3}{2}\Delta+1 \big\rceil$, the Markov chain $\C(\CH)$ has
mixing time $O(n\log n)$.
\end{thm}
\begin{proof}
Consider two proper colourings $X$ and $Y$ differing in a single
vertex $w$. Without loss of generality let the change vertex $w$ be
coloured 1 in $X$ and 2 in $Y$. We will partition the edges $e\in
\mathcal{E}$ containing $w$ into four classes $E_1, E_2, E_3, E_4$,
determined by the colouring of $e\setminus \{w\}$, as follows:
\[E_1: \{1,2\}\qquad  E_2: \{1,i\}\ \mathrm{or}\ \{2,i\}\ \
(2<i)\qquad E_3: \{i,i\}\ \ (2<i)\qquad  E_4: \{i,j\}\ \ (2<i< j).\]
Instead of using Hamming distance, we will take a new metric defined
by
\[ \dd(X,Y)=\sum_{i=1}^{4}c_{i} |E_i|,\]
where $1=c_1\geq c_2\geq c_3\geq c_4>0$, and for convenience
$c_0=\Delta+1$. Note that $\dd(X,Y)\leq\Delta$ if $X,Y$ have Hamming
distance 1. The diameter is therefore at most $\Delta n$ in the
metric $\dd$.

Arguing as in Section~\ref{sec:indsets}, we have
\renewcommand{\arraystretch}{1.3}
\begin{equation}\label{col-eq10}
\begin{array}{rl}
nq \E[\dd_1-\dd_0]\ \leq & -(q-|E_3|)\big(c_1|E_1|+c_2|E_2|+c_3|E_3|+c_4|E_4|\big)\\
& + |E_1|\big(-2(q-\Delta-1)(c_{1}-c_2)+2(c_{0}-c_1)\big)\\
& + |E_2|\big(-(q-\Delta-2)(c_{2}-c_4)-(c_{2}-c_3)+(c_{0}-c_2)+(c_{1}-c_2)\big)\\
& + |E_3|\big(-2(q-\Delta-2)(c_{3}-c_4)+4(c_{2}-c_3)\big)\\
& + |E_4|\big(2(c_{3}-c_4)+4(c_{2}-c_4)\big).
\end{array}
\end{equation}
If, in (\ref{col-eq10}), we set
\begin{equation}\label{col-eq20}
\begin{array}{rl}
2(q-\Delta-1)(c_{1}-c_2)-2(c_{0}-c_1)+c_1(q-|E_3|)&\ =\ \gamma\\
(q-\Delta-2)(c_{2}-c_4)+(c_{2}-c_3)-(c_{0}-c_2)-(c_{1}-c_2)+c_2(q-|E_3|)&\ =\ \gamma\\
2(q-\Delta-2)(c_{3}-c_4)-4(c_{2}-c_3)+c_3(q-|E_3|)&\ =\ \gamma\\
-2(c_{3}-c_4)-4(c_{2}-c_4)+c_4(q-|E_3|)&\ =\ \gamma,
\end{array}
\end{equation}
where $\gamma\geq 0$, we have
\begin{equation}\label{col-eq30}
 \E[\dd_1] \leq\  \dd_0-\frac{\gamma\Delta}{nq} \leq \Big(1-\frac{\gamma}{nq}\Big)\dd_0.
\end{equation}
Note, that if we put $q'=q-|E_3|$, $\Delta'=\Delta-|E_3|$ in
(\ref{col-eq20}), we have
\begin{equation}\label{col-eq40}
\begin{array}{rl}
2(q'-\Delta'-1)(c_{1}-c_2)-2(c_{0}-c_1)+c_1q'&\ =\ \gamma\\
(q'-\Delta'-2)(c_{2}-c_4)+(c_{2}-c_3)-(c_{0}-c_2)-(c_{1}-c_2)+c_2q'&\ =\ \gamma\\
2(q'-\Delta'-2)(c_{3}-c_4)-4(c_{2}-c_3)+c_3q'&\ =\ \gamma\\
-2(c_{3}-c_4)-4(c_{2}-c_4)+c_4q'&\ =\ \gamma.
\end{array}
\end{equation}
This corresponds to a system like (\ref{col-eq20}) with degree
$\Delta'$, $q'$ colours and $|E_3|=0$. But, since
$q'/\Delta'=(q-|E_3|)/(\Delta-|E_3|)\geq q/\Delta$, the smallest
ratio for $q/\Delta$ is given by setting $|E_3|=0$ in
(\ref{col-eq20}). Also, putting $c_{3}=c_4$ makes the third and
fourth equations in (\ref{col-eq40}) identical, so $c_{3}=c_4$ must
be a solution. With these simplifications, and putting
$c_{0}=\Delta+1$, $c_{1}=1$, we have
\begin{equation*}
\begin{array}{rl}
2(q-\Delta-1)(1-c_2)-2\Delta+q&\ =\ \gamma\\
(q-\Delta-1)(c_{2}-c_4)-2(1-c_2)-\Delta+c_2q&\ =\ \gamma\\
-4(c_{2}-c_4)+c_4q&\ =\ \gamma.
\end{array}
\end{equation*}
Now the linear equations (\ref{col-eq20}) may be  solved for
$c_{2}$, $c_4$ and $\gamma$, giving
\begin{equation*}
c_{1}=1,\quad c_{2}=\frac{2q-2\Delta+1}{2q-\Delta+1},\quad
c_{3}=c_4=\frac{2q-3\Delta+1}{2q-\Delta+1},\quad
\gamma=\frac{2q^{2}-q(3\Delta-1)-4\Delta}{2q-\Delta+1}.
\end{equation*}
The condition $\gamma\geq 0$ is equivalent to
\begin{equation*}
   q\ \geq\
   \tfrac{3\Delta-1}{4}\Big(1+\sqrt{1+\tfrac{32\Delta}{(3\Delta-1)^2}}\Big),\qquad
   \mathrm{i.e.}\ \ q\ \geq\\
     \big\lceil \tfrac{3}{2}\Delta\big\rceil+1 .
\end{equation*}
Note that we have $c_{i}>0$ ($i=1,\ldots,4$) under this condition.
Note also that $\gamma>0$ and hence, using (\ref{col-eq30}), the
mixing time satisfies
\[\tau(\varepsilon)\ \leq\ \frac{2q^{2}-q\Delta+q}{2q^{2}-q(3\Delta-1)-4\Delta}
\ n \ln \Big(\frac{\Delta n}{\varepsilon}\Big).\qedhere\]
\end{proof}

\section{Colouring bipartite graphs}\label{sec:bipartite}
Our final application is to colouring bipartite graphs. Several
recent papers have used a stopping times or ``burn in'' analysis to
prove rapid mixing for Glauber dynamics of graph colouring, starting
with~\cite{DF03}. These are largely based upon the idea that
although a vertex can have only $q-\Delta$ colours with which
to be properly recoloured, it is very unlikely for any vertex to have so
few colours available after a period of ``burn in''. Subject to more stringent girth and degree restrictions than used here, rapid mixing has been proved for fewer colours~\cite{DFHV04, HV03, Mol02}. Here we capture
this intuition by using a metric which
directly incorporates the number of colours available to a vertex.
In order to simplify the analysis, we do not consider Glauber
dynamics here. Instead we prove that a Markov chain \textsc{Scan} which uses the same method for recolouring a vertex as Glauber dynamics, but recolours the vertices in a deterministic order, mixes rapidly. In order to show this we first prove results for a closely related Markov chain, \textsc{Multicolour}, which is of interest in its own right.

Let $G=(V, E)$ be a bipartite graph with bipartition $V_1,V_2$, and
maximum degree $\Delta$. For $v \in V$, let $\mathcal{N}(v)=\{w:
\{v,w\}\in E\}$ denote the neighbourhood of $v$. Let $Q=[q]$ be a colour set,
and $X: V\rightarrow Q$ be a colouring of $G$, not necessarily
proper. Let $C(v)=\{ X(w):w\in \mathcal{N}(v)\}$ be the set of colours occurring in the neighbourhood of $v$, and $c(v)=|C(v)|$. We consider the Markov chain \textsc{Multicolour} on
colourings of $G$, which in each step picks one side of the
bipartition at random, and then recolours every vertex on that side, followed by recolouring every vertex in the other half of the bipartition.
If the state of \textsc{Multicolour} at time $t$ is $X_t$, the state
at time $t+1$ is given by\vspace{2ex}

\begin{minipage}[l]{0.9\textwidth}
\textsc{Multicolour}\vspace{-1ex}
\begin{enumerate}
  \item choosing  $r \in \{1,2\}$ uniformly at random,\label{item1}\vspace{-1ex}
  \item for each vertex $v \in V_r$, \label{item2}\vspace{-1ex}
\begin{enumerate}
  \item choosing a colour $q(v)\in Q\backslash C(v)$ uniformly at random, 
  \item setting $X_{t+1}(v)= q(v)$. (Heat bath recolouring)\label{item3}
\end{enumerate}
  \item for each vertex $v \in V \backslash V_r$, \label{item4}\vspace{-1ex}
\begin{enumerate}
  \item choosing a colour $q(v)\in Q\backslash C(v)$ uniformly at random, 
  \item setting $X_{t+1}(v)= q(v)$. 
\end{enumerate}
\end{enumerate}
\end{minipage}\vspace{1ex}

Note that the order in which the vertices are processed in
steps~\ref{item2} and~\ref{item4} is immaterial. This chain is a single-site dynamics intermediate between Glauber
and scan. It is easy to see that it is ergodic if $q>\Delta+1$, and
has equilibrium distribution uniform on all proper colourings of
$G$. Observe also that it requires considerably fewer random bits
than Glauber, and only slightly more than scan. We prove the following theorem.
\begin{thm}\label{bipcol}
For $q>f(\Delta)$ the mixing times of \textsc{Scan} and \textsc{Multicolour} are $O(\log(n))$, where $f$ is a function such that
\begin{enumerate}
\item $f(\Delta)\rightarrow \beta \Delta$, as $\Delta\rightarrow \infty$, where $\beta$ satisfies $\tfrac{1}{\beta} e^{1/\beta}=1$,
\item $f(\Delta)\leq \lceil 11 \Delta/6 \rceil$ for $\Delta \geq 14$,
\item $f(\Delta) < \lceil 11 \Delta/6 \rceil$ for $\Delta \geq 31$,
\item in particular $f(22) =40 < \lceil 11 \Delta/6 \rceil$.
\end{enumerate}
\end{thm}
We will require the following lemmas.

\begin{lem}\label{dflem}
For $1\leq i \leq \Delta$ let $S_i$ be a subset of $(Q-q_0)$ such that $m_i=|S_i|\geq q-\Delta$. Let $s_i$ be selected uniformly at random from $S_i$, independently for each $i$. Finally let $C=\{s_i:1\leq i \leq \Delta\}$ and $c=|C|$. Then \[\E[q-c\ |\ s_1=q_1]\geq 1+ (q-2)\left(   1-\frac{1}{q-\Delta}\right)^{\frac{(\Delta-1)(q-\Delta)}{q-2}} =\alpha.\]
\end{lem}
\begin{proof}
 This follows from~\cite[Lemma 2.1]{DF03} with minor adjustments as follows. Let $a_{ij}=1$ if $j\in S_i$ and $0$ otherwise. Thus $m_i=\sum_{j\in (Q-q_0)} a_{ij}$ and 
\[
\E[q-c] = 1+ \sum_{j\in (Q-q_0)} \prod_{i=1}^\Delta \left(1-\frac{1}{m_i}\right)^{a_{ij}}.
\]
However if we are given that $s_1=q_1$, then 
\begin{align*}
\E[q-c\ |\ s_1=q_1] &=\ 1+ \sum_{j\in (Q-q_0-q_1)} \prod_{i=2}^\Delta \left(1-\frac{1}{m_i}\right)^{a_{ij}}\\
&\geq\ 1+ (q-2)\left(  \prod_{j\in (Q-q_0-q_1)} \prod_{i=2}^\Delta \left(1-\frac{1}{m_i}\right)^{a_{ij}} \right)^{\frac{1}{q-2}}\\
&\geq\ 1+ (q-2)\left(  \prod_{i=2}^\Delta \left(1-\frac{1}{m_i}\right)^{m_i} \right)^{\frac{1}{q-2}}\\
&\geq\ 1+ (q-2)\left(   1-\frac{1}{q-\Delta}\right)^{\frac{(\Delta-1)(q-\Delta)}{q-2}}. 
\end{align*}
Where the final inequality follows because $(1-1/m_i)^{m_i}$ in increasing with $m_i$ and $m_i\geq q-\Delta$ for all $i$. 
\end{proof}

\begin{lem}\label{dflem2}
For $1\leq i \leq \Delta$ let $S_i$ be a subset of $(Q-q_0)$ such that $m_i=|S_i|\geq q-\Delta$. Let $s_i$ be selected uniformly at random from $S_i$, independently for each $i$. Finally let $C=\{s_i:1\leq i \leq \Delta\}$ and $c=|C|$. Then
\[\E\left[\frac{1}{q-c} \ | \  s_1=q_1 \right]\leq \frac{1}{\alpha}\left(1+\frac{(q-\alpha-1)(\alpha-1)}{(q-\Delta)(q-2)\alpha}\right)=\alpha'.\]
\end{lem}
 \begin{proof}
We will write $\bar{c}$ for $\E[c\ | \  s_1=q_1]$. Let $Z=\frac{c-\bar{c}}{q-c}$, so that 
\begin{equation}\label{q-c1}
\frac{1}{q-c}=\frac{1}{q-\bar{c}}\left(\frac{1}{1-Z}\right).
\end{equation}
 Note that $(1-Z)^{-1}=\tfrac{q-\bar{c}}{q-c}\leq \tfrac{q-\bar{c}}{q-\Delta}$. Now \[
\frac{1}{1-Z}= 1+Z+\frac{Z^2}{1-Z}\leq 1+Z+\frac{(q-\bar{c})Z^2}{q-\Delta}.\] 
Hence 
\begin{equation}\label{q-c2}
\E[(1-Z)^{-1}\ | \  s_1=q_1]\leq 1+\frac{q-\bar{c}}{q-\Delta}\frac {\textrm{Var} (c\ | \  s_1=q_1)}{(q-\bar{c})^2}=1+\frac {\textrm{Var} (c\ | \  s_1=q_1)}{(q-\Delta)(q-\bar{c})}.\end{equation}

We now turn our attention to bounding $\textrm{Var} (c\ | \  s_1=q_1)$. Let $c=\sum_{j\in (Q-q_0)} I_j$, where $I_j$ indicates that colour $j$ is in $C$. Now, conditional on $ s_1=q_1$, we have  \[\textrm{Var} (\sum_{j\in (Q-q_0)} I_j)=\sum_{j\in (Q-q_0)}\textrm{Var} ( I_j)+ 2 \sum_{j<k}\textrm{Cov}(I_j,I_k)\leq \sum_{j\in (Q-q_0)}\textrm{Var} ( I_j),\]
since $I_j$ and $I_k$ are negatively correlated for all $j$ and $k$. Let $p_j=\Pr(I_j=1)$, then $I_j$ has variance $p_j(1-p_j)$ and $\sum_{j\in (Q-q_0)} p_j=\bar c$. Also note that $p_{q_1}=1$, hence $\textrm{Var}(I_{q_1})=0$. By convexity, the maximum of $\sum_{j\in (Q-q_0-q_1)} p_j(1-p_j)$ such that $\sum_{j\in (Q-q_0-q_1)} p_j=\bar c -1$ is given by setting $p_j=(\bar{c}-1)/(q-2)$. Hence, using $\bar c= q-\alpha$, 
\begin{equation}\label{q-c3}
\textrm{Var} (c\ | \  s_1=q_1)\leq (\bar{c}-1)\left(1-\frac{\bar{c}-1}{(q-2)}\right)=\frac{(q-\alpha-1)(\alpha-1)}{q-2}.\end{equation}

Putting together equations~(\ref{q-c1}), (\ref{q-c2}) and (\ref{q-c3}) we have \[\E\left[\frac{1}{q-c}\ | \  s_1=q_1\right]\leq \frac{1}{\alpha}\left(1+\frac{(q-\alpha-1)(\alpha-1)}{(q-2)}\frac{1}{(q-\Delta)\alpha}\right).\]
 \end{proof}

\begin{proof}[Proof of Theorem~\ref{bipcol}] We first prove the theorem for \textsc{Multicolour}. In the path coupling setting, we will take $S$ to be the set of
pairs colourings which differ at exactly one vertex.  Let $v$ be the
change vertex for some pair $(X,Y)\in S$, and assume without loss
that $v\in V_1$. The distance between $X$ and $Y$ is defined to be $\dd(X,Y)=\sum_{w\in \mathcal{N}(v)} \frac{1}{q-c_{X,Y}(w)}$, where $c_{X,Y}(w)$ is taken to be $\min\{c_X(w),c_Y(w)\}$ in the case that they differ. We couple as follows (the usual path coupling for
Glauber dynamics). If we are recolouring a vertex which is not a neighbour of $v$, then the sets of available colours in $X$ and $Y$ are the same, and we use the same colour in both copies of the chain. If we are recolouring a vertex $w\in \mathcal{N}(v)$ then there are three cases to consider:
\begin{enumerate}
\item \label{RBblocked}$|\{X(v),Y(v)\}\cap \{X(z):z\in \mathcal{N}(w)\backslash \{v\}\}|=2.$ \\ The colours $X(v)$ and $Y(v)$ are not available for recolouring $w$ in either copy of the chain, hence the sets of available colours are the same, and we use the same colour in both copies of the chain.
\item $|\{X(v),Y(v)\}\cap \{X(z):z\in \mathcal{N}(w)\backslash \{v\}\}|=1.$ \\ Without loss assume colour $X(v)$ is not available to $w$ in either copy of the chain. Colour $Y(v)$ is only available in $X$. We couple recolouring $w$ in $X$ with any colour other than $Y(v)$, with recolouring using the same colour in $Y$. We couple recolouring $w$ in $X$ with colour $Y(v)$, uniformly between recolouring $w$ with each available colour in $Y$.
\item $|\{X(v),Y(v)\}\cap \{X(z):z\in \mathcal{N}(w)\backslash \{v\}\}|=0.$ \\ Here colour $Y(v)$ is only available in chain $X$, and $X(v)$ in only available in $Y$. We couple together recolouring with these colours respectively, and for each other colour (that is available to both copies), we recolour $w$ with the same colour in both $X$ and $Y$.
\end{enumerate}
Note that in case \ref{RBblocked}, there is no probability of $w$ being coloured differently in the two chains. In the other cases, the probability of disagreement at $w$ is $\frac{1}{q-c_{X,Y}(w)}$.

Let $X',Y'$ be the colourings after recolouring $V_r$ (half a step of \textsc{Multicolour}) and $X'',Y''$ be the colourings after the full step of  \textsc{Multicolour}. If we randomly select $V_1$ to be recoloured first, then the two copies of the chain have coupled in $X'$ and $Y'$ since the vertices in $V_1$ have the same set of available colours in each chain. 

So suppose that we select $V_2$ to be recoloured first. The only vertices in $V_2$ that have different sets of available colours are those which are neighbours of $v$. Let $\mathcal N(v)= \{w_1,\ldots,w_k\}$ and consider the path $W_0, W_1,\ldots, W_{k+1}$ from $X'$ to $Y'$, where for $1\leq i\leq k$, $W_i$ agrees with $X'$ on all vertices except $w_1,\ldots, w_i$ which are coloured as in $Y'$, and $W_0=X'$ and $W_{k+1}=Y'$. Then for $i\leq k$ we have \begin{equation} \label{eqnwi-1wi}
\dd(W_{i-1},W_i)=\mathds{1}_{w_i}\sum_{z\in \mathcal{N}(w_i)} \frac{1}{q-c_{W_{i-1},W_i}(z)}\leq \mathds{1}_{w_i}\sum_{z\in \mathcal{N}(w_i)} \frac{1}{q-c_{W_i}(z)},\end{equation}
 where $\mathds{1}_{w_i}$ indicates whether $X'$ and $Y'$ differ on $w_i$. Note that $\Pr[\mathds{1}_{w_i}=1]\leq \frac{1}{q-c_{X,Y}(w_i)}$. Furthermore, by the construction of the coupling either conditioning on $\mathds{1}_{w_i}=1$ is the same as conditioning that $W_{i-1}(w_i)=q_1$, or that $W_{i}(w_i)=q_1$, for some $q_1$. We assume without loss that this is $W_i$. Then for each $z\in \mathcal N(w_i)-v$ the selection of colours in $C_{W_i}(z)$ satisfies the conditions of Lemma~\ref{dflem2}, since we may take $q_0=X(z)$ and $q_1$ as above.  For $v$, there is no colour $q_0$ which is necessarily unavailable for all its neighbours, since some are coloured as in $X'$ and some as in $Y'$. Hence we use a slightly weaker bound on $\alpha$ and $\alpha'$, given by \[
 \alpha_v=(q-1)\left(   1-\frac{1}{q-\Delta}\right)^{\frac{(\Delta-1)(q-\Delta)}{q-1}}\ \textrm{ and }\quad
 \alpha'_v=\frac{1}{\alpha_v}\left(1+\frac{(q-\alpha_v)(\alpha_v)}{(q-\Delta)(q-1)\alpha_v}\right).
 \] Hence for $i\leq k$,  $\E[\dd(W_{i-1},W_i)]\leq  \frac{1}{q-c_{X,Y}(w_i)}((\Delta-1)\alpha'+\alpha'_v)$. The value of $\dd(W_k,W_{k+1})$ is still $\dd(X,Y)$ since the vertices in $V_1$ have not yet been recoloured.

\begin{table}[htbp]
   \centering
   \begin{tabular}{@{} lccr @{}} 
      \toprule
      $\Delta$    & $q$ &$\lceil 11\Delta/6  \rceil$ & $q/\Delta$\\
      \midrule
9 & 17 & 17 & 1.89 \\
10 & 19 & 19 & 1.90 \\
11 & 21 & 21 & 1.91 \\
12 & 23 & 22 & 1.92 \\
13 & 25 & 24 & 1.92 \\
14 & 26 & 26 & 1.86 \\
15 & 28 & 28 & 1.87 \\
16 & 30 & 30 & 1.88 \\
17 & 32 & 32 & 1.88 \\
18 & 33 & 33 & 1.83 \\
19 & 35 & 35 & 1.84 \\
20 & 37 & 37 & 1.85 \\
21 & 39 & 39 & 1.86 \\
22 & 40 & 41 & 1.82 \\
23 & 42 & 43 & 1.83 \\
24 & 44 & 44 & 1.83 \\
25 & 46 & 46 & 1.84 \\
26 & 48 & 48 & 1.85 \\
27 & 49 & 50 & 1.81 \\
28 & 51 & 52 & 1.82 \\
29 & 53 & 54 & 1.83 \\
30 & 55 & 55 & 1.83 \\
31 & 56 & 57 & 1.81 \\
32 & 58 & 59 & 1.81 \\
33 & 60 & 61 & 1.82 \\
34 & 61 & 63 & 1.79 \\
35 & 63 & 65 & 1.80 \\
36 & 65 & 66 & 1.81 \\
37 & 67 & 68 & 1.81 \\
38 & 68 & 70 & 1.79 \\
39 & 70 & 72 & 1.79 \\
40 & 72 & 74 & 1.80 \\
41 & 74 & 76 & 1.80 \\
42 & 75 & 77 & 1.79 \\
43 & 77 & 79 & 1.79 \\
44 & 79 & 81 & 1.80 \\
45 & 81 & 83 & 1.80 \\
46 & 83 & 85 & 1.80 \\
47 & 84 & 87 & 1.79 \\
48 & 86 & 88 & 1.79 \\
49 & 88 & 90 & 1.80 \\
50 & 90 & 92 & 1.80 \\
10000 & 17634 & 18334 & 1.76 \\
      \bottomrule
   \end{tabular}
   \caption{Minimum values of $q$ for contraction.}
   \label{table}
\end{table}

Now we consider the vertices in $V_1$.  We apply the same analysis as above to each path segment $W_{i-1},W_i$, but augment the analysis using the fact that at the time a vertex $z\in V_1$ is recoloured, its neighbours (in $V_2$) will already have been randomly recoloured. Let the neighbours of $w_i$ be $z_1, z_2, \ldots z_l$, and consider the path $Z_0, Z_1, \ldots Z_{l+1}$, where for $1\leq j\leq l$, $Z_{j}$ agrees with $W_{i-1}$ on all vertices except $z_1,\ldots, z_j$ which are coloured as in $W_i$, and $Z_0=W_{i-1}$ and $Z_{l+1}=W_i$. Arguing as above, for $j\leq l$ we have  
\[
\dd(Z_{j-1},Z_{j})=\mathds{1}_{z_j}\sum_{w\in \mathcal{N}(z_j)} \frac{1}{q-c_{Z_{i-1},Z_i}(w)}. 
\]
But now $\Pr[\mathds{1}_{z_j}=1 |\  W_{i-1},W_i]\leq \frac{1}{q-c_{W_{i-1},W_i}(z_j)}\mathds{1}_{w_i} $. This is similar to equation~(\ref{eqnwi-1wi}), and the same argument gives $\E[\mathds{1}_{z_j}=1]\leq \frac{1}{q-c_{X,Y}(w_i)}\alpha'$, for $z_j\neq v$ and $\E[\mathds{1}_{z_j}=1]\leq \frac{1}{q-c_{X,Y}(w_i)}\alpha'_v$ if $z_j=v$. 
Also, since it depends only on the colouring of $V_2$, we have $\dd(Z_l,Z_{l+1})=\dd(W_{i-1},W_i)$. So 
\[
\E[\sum_{j=1}^{l+1} \dd(Z_{j-1},Z_j)] \leq \frac{1}{q-c_{X,Y}(w_i)}((\Delta-1)\alpha'+\alpha'_v)(((\Delta-1)\alpha'+\alpha'_v)+1).
\]
Finally note that $W_k$ and $W_{k+1}$ differ only in $V_1$, so after recolouring $V_1$ they have coupled. Hence 
\begin{align}
\E[\dd(X'',Y'')]&=\ \frac12 \sum_{i=1}^k \sum_{j=1}^{l+1}\E[\dd(Z_{j-1},Z_j) ]\\
&\leq\ \frac12 \sum_{i=1}^k \frac{1}{q-c_{X,Y}(w_i)}((\Delta-1)\alpha'+\alpha'_v)( ((\Delta-1)\alpha'+\alpha'_v)+1)\\
&=\ \dd(X,Y)((\Delta-1)\alpha'+\alpha'_v) \frac{( ((\Delta-1)\alpha'+\alpha'_v)+1)}{2}.
\end{align}
This gives contraction as long as $((\Delta-1)\alpha'+\alpha'_v)$ is less than $1$. For large $\Delta$, we see that $\alpha'$ and $\alpha'_v$ both approach $\tfrac1q e^{\Delta/q}$. Hence we have contraction when $\tfrac\Delta q e^{\Delta/q}<1$. For small values of $\Delta$ it is possible to compute the smallest integral value of $q$ for which there is contraction. These values are shown in Table~\ref{table}. When there is contraction, standard path coupling arguments give the mixing time bounds claimed.

We now argue that \textsc{Scan} mixes as rapidly as \textsc{Multicolour}. The Markov chain \textsc{Scan} recolours the two sides of the bipartition in order, $(V_1,V_2),(V_1,V_2)\ldots$. The Markov chain \textsc{Multicolour} recolours a random side first in each step. However, recolouring the same side twice in a row has exactly the same effect as recolouring it once, since vertices in the same side of the bipartition are independent. The recolouring given by a run of \textsc{multicolour} with order $(V_1,V_2), (V_2,V_1),(V_1,V_2)$ has exactly the same result as if the reversed pair was omitted. Hence any randomly chosen sequence can be replaced with a purely alternating sequence. Should the purely alternating sequence corresponding to the random choices of \textsc{Multicolour} start with $V_2$ or finish with $V_1$, we can augment the sequence with a recolouring of $V_1$ at the beginning or $V_2$ at the end respectively. The result follows, since the former is equivalent to taking a different starting position in \textsc{Multicolour}, and the latter cannot increase the total variation distance from stationarity.
\end{proof}

\begin{rem}
Our analysis shows that one-step analysis of a single-site chain on
graph colourings need not break down at $q=2\Delta$~\cite{J95,SS97}.
This apparent ``boundary'' seems merely to be an artefact of using
Hamming distance.
\end{rem}
\begin{rem}
Our scan chain can be used to prove polynomial mixing time for the
Glauber dynamics (with the same values of $q$ and $\Delta$) by comparison
techniques~\cite{DS93,RT00}. However,
the proof is not completely straightforward and will appear elsewhere. 
\end{rem}
\begin{rem}
We note that many of the infinite graphs studied in statistical physics are
bipartite, for example cubic grids and trees.
Therefore our results imply, for example, absence of phase transition in 
the antiferromagnetic Potts model in the cubic grid with $q$ colours
and dimension $d=\Delta/2$. A proof follows the
lines of that given by Vigoda~\cite[\S 5]{V99} with obvious modifications.
Since results with similar $q,d$ have been proved by different arguments
in~\cite{GMP04}, we omit the details.
\end{rem}

\end{document}